\documentclass{article}

\usepackage{PRIMEarxiv}

\usepackage[utf8]{inputenc} % allow utf-8 input
\usepackage[T1]{fontenc}    % use 8-bit T1 fonts
\usepackage{hyperref}       % hyperlinks
\usepackage{url}            % simple URL typesetting
\usepackage{booktabs}       % professional-quality tables
\usepackage{amsfonts}       % blackboard math symbols
\usepackage{nicefrac}       % compact symbols for 1/2, etc.
\usepackage{microtype}      % microtypography
\usepackage{lipsum}
\usepackage{fancyhdr}       % header
\usepackage{graphicx}       % graphics
\graphicspath{{media/}}     % organize your images and other figures under media/ folder
\usepackage{graphicx}
\usepackage{amsmath}
\usepackage{amsfonts}
\usepackage{amssymb}
\usepackage{times}
\usepackage{braket}
\usepackage{url}

\def\boldtheta#1{\boldsymbol{\theta}^{#1}}

\def\boldx#1{\boldsymbol{x}#1}
\def\tildey{\widetilde{\boldsymbol{y}}}
\def\tildePhi{\widetilde{\boldsymbol{\Phi}}}
\def\bolde#1{\boldsymbol{e}_{#1}}
\def\hatRho{\widehat{{\rho}}}
\def\bigSigma{\text{\large $\Sigma$}}

\title{Input-output reduced order modeling for public health intervention evaluation
}

\author{
  Alex Viguerie$^1$, Md Hafizul Islam$^3$, Evin Uzun Jacobson$^4$ \\
  Division of HIV Prevention \\
  Centers for Disease Control and Prevention \\
  Atlanta, GA (USA), 30333\\
  \texttt{$^1$xjm9@cdc.gov, $^3$qla0@cdc.gov, $^4$wqm4@cdc.gov } \\
   \And
  Chiara Piazzola$^2$ \\
  Department of Mathematics \\
  Technical University of Munich \\
  85748, Garching bei M{\"u}nchen, Germany\\
  \texttt{chiara.piazzola@tum.de} \\
}

\begin{document}
\maketitle

\begin{abstract}
In recent years, mathematical models have become an indispensable tool in the planning, evaluation, and implementation of public health interventions. Models must often provide detailed information for many levels of population stratification. Such detail comes at a price: in addition to the computational costs, the number of considered input parameters can be large, making effective study design difficult. To address these difficulties, we propose a novel technique to reduce the dimension of the model input space to simplify model-informed intervention planning. The method works by first applying a dimension reduction technique on the model output space. We then develop a method which allows us to map each reduced output to a corresponding vector in the input space, thereby reducing its dimension. We apply the method to the HIV Optimization and Prevention Economics (HOPE) model, to validate the approach and establish proof of concept.
\end{abstract}

% keywords can be removed

\section{Introduction}
In recent years, mathematical models have become an indispensable tool in the planning, evaluation, and implementation of public health interventions. Models often need to provide detailed information for many distinct levels of population stratification. Such detail comes at a price: in addition to computational costs, the number of considered input parameters can be large. This leads to difficulty at all stages of the modeling workflow. For example, effective study designs can be difficult to develop and computationally demanding, and optimization algorithms are slow to converge and may exhibit (perhaps infinitely) many local solutions. Post-processing of model results can also be challenging, as the large number of varied parameters can make drawing concrete, actionable conclusions from simulation results difficult. 

Recent years have seen significant progress in \textit{reduced-order modeling} (ROM) techniques to reduce the dimension of mathematical models. Particularly, popular techniques include approaches based on the Proper Orthogonal Decomposition (POD) \cite{luo2018proper,quarteroni2014reduced}, surrogate modeling \cite{jakeman2022adaptive}, and Dynamic Mode Decomposition (DMD) and related methods based on Koopman theory \cite{baddoo2023physics,hess2023data}.

Such methods have been applied in public health for short-term forecasting and diagnostic applications\cite{viguerie2022coupled}. However, intervention planning generally requires multiple model evaluations for different input configurations. Hence, ROM techniques that are useful for forecasting purposes, may not be applicable. Instead, intervention planning involves observing how model outputs \textit{change} in response to variation in model inputs. Thus, using ROMs for intervention planning requires an approach that can incorporate parameter-dependence. While ROM techniques have been developed for parameterized problems, many of these applications are focused on cases in which the primary difficulty arises from the high dimension of the \textit{output} space\cite{luo2018proper,quarteroni2014reduced,hess2023data,andreuzzi2023dynamic}. 

However, in intervention planning, most of the difficulty arises from the dimension of the \textit{input space}, a problem which has received less attention. A recent work \cite{kim2024dimensionality} explored the issue of extracting surrogate models from problems with high-dimensional input spaces, using a combination of dimension reduction techniques and Bayesian distribution fitting.  Other notable works in this area focused on constructing surrogate models \cite{gadd2019surrogate}, and input-output ROM schemes for experimental data \cite{ma2011kernel}. Active subspace (AS) and related methods have also been used to reduce input dimensionality for scalar-valued problems \cite{constantine2015active}, and have been applied to problems in public health\cite{romor2022kernel}. However, AS methods are applicable to scalar-valued functions, and applications to vector-valued problems typically require the use of a scalar-valued surrogate; for example, using the basic reproduction number as a surrogate for a model of infectious disease transmission\cite{romor2022kernel}.

We propose herein a novel method for input parameter-space reduction. The method works by applying dimension-reduction techniques over an aggregation of model outputs, generated over a range of input parameters. Through surrogate modeling over the reduced-order space, we construct a mapping which maps each reduced-order output to a vector of input parameters in a reduced-order input space. The reduced-order inputs can then be used directly to identify intervention plans based on program targets and constraints.

The paper is outlined as follows. In Section 2, we briefly introduce our model problem. In Section 3, we explain the proposed method for input/output model order reduction. In Section 4, we demonstrate the potential utility of the method by applying it on a large-scale model of HIV transmission, the HIV Optimization and Prevention Economics (HOPE) model. We conclude by summarizing our findings and discussing directions for future work.

\section{Basic model and preliminaries}
\label{sec:headings}

Let us consider a nonlinear continuous-time dynamical system consisting of $n$ unknowns and depending on an $m$-dimensional vector $\widehat{\boldtheta{}} = [\widehat{\theta}_1, \ldots, \widehat{\theta}_m]\in \mathbb{R}^m$ of parameters, assumed constant in time. Denote as $\Gamma=\Gamma_1 \times \Gamma_2 \times ... \times \Gamma_m$ the set of all possible values of $\widehat{\boldtheta{}}$. Each $\Gamma_i$ is assumed to be closed and bounded herein, though this assumption is not strictly necessary.
\par We may then write:
\begin{align}\begin{split}\label{genericNonlinearSystem}
\dfrac{d\boldx{}}{dt} &= F(t,\boldx{},\widehat{\boldtheta{}}) \qquad \text{ for }  t_0< t \leq t_{end}, \\
\boldx{} &= \boldx{}_0 \qquad\qquad\qquad \text{ at } t=t_0,
\end{split}\end{align}
with $F$ a nonlinear function depending on $\widehat{\boldtheta{}}$, as well the solution state $\boldx$ and time $t$. We assume that the system \eqref{genericNonlinearSystem} is sufficiently regular such that the solution $\boldx{}$ is time-continuous for all $\widehat{\boldtheta{}}\in \Gamma$ and that the solution is continuous in $\widehat{\boldtheta{}}$ at all $t$ for all $\widehat{\boldtheta{}} \in \Gamma$.
\par Let us assume that $\widehat{\theta}_i$, $i=1,\ldots,m$ are independently distributed random variables taking values in $\Gamma_i$ and  $\widehat{\rho}_i(\widehat{\theta}_i):\Gamma_i \to \lbrack 0,\,\infty)$ be the corresponding probability density functions. Then the joint PDF $\hatRho: \Gamma \rightarrow [0,\infty)$ of $\widehat{\boldtheta{}}$ is given by:
\begin{align}\label{jointPDFBoldTheta}
\hatRho(\widehat{\boldtheta{}}) &:= \prod_{i=1}^{m} \widehat{\rho}_i(\widehat{\theta}_i).
\end{align}
For reasons that will be clear later, we want to consider all of the input parameters as defined over a uniform interval. Let $\boldsymbol{\varphi} : \Gamma\to \lbrack 0,\,1\rbrack^m $ be a bijective map. Hence each $\widehat{\boldtheta{}} \in \Gamma $ can be associated with a unique $\boldtheta{}\in \lbrack0,\,1\rbrack^m$ such that:
\begin{align}\label{bijectiveMap}
    \widehat{\boldtheta{}} = \boldsymbol{\varphi}^{-1}(\boldtheta{}).
\end{align}
Combining \eqref{jointPDFBoldTheta} and \eqref{bijectiveMap}, the joint PDF $\rho$ of $\boldtheta{}$ is given by:
\begin{equation}\label{jointPDFTheta}
\rho(\boldtheta{}) = \left(\hatRho \circ \boldsymbol{\varphi}^{-1}\right)(\boldtheta{}) =  \prod_{i=1}^{m} \left(\widehat{\rho}_i \circ \varphi_i^{-1}  \right)({\theta}_i).
\end{equation}
Let $\{ \boldtheta{j}\}_{j=1}^k $ be $k$ samples of $\boldtheta{}$. Define the \textit{vector-valued function} $\boldsymbol{y}(\boldtheta{}):\lbrack 0,\,1\rbrack^m \to \mathbb{R}^d$ to account for outcomes of interest of \eqref{genericNonlinearSystem}, whose $i$-th entry is given by:
\begin{align}\label{defineFunctional}
    y_i (\boldtheta{}) = J_i \left(\boldx{}\circ \boldsymbol{\varphi}^{-1}\right) (\boldtheta{}),
\end{align}
 with $J_i$ a bounded and continuous functional, possibly nonlinear. For example, the $J_i$ may represent the annual disease incidence or disease-related mortality among an intervention population in a given time period. Note that $\boldsymbol{y}(\boldtheta{})$ does not depend on time $t$ directly. One may introduce time-dependent quantities through the definition of $J_i$, for example two distinct $J_i$ may describe the same output considered over different time periods.

\section{Input-output model order reduction}
\subsection{Reduction of the output space}
Given $k$ samples of $\boldtheta{}$, we solve \eqref{genericNonlinearSystem} for each $\widehat{\boldtheta{}}^j=\boldsymbol{\varphi}^{-1}(\boldtheta{j})$. Denote $\boldsymbol{y}^j= \boldsymbol{y}(\boldtheta{j})$. Let $Y$ be the $d \times k$ \textit{snapshot matrix}:
\begin{align}\label{snapshotMtx}
Y= \begin{bmatrix} \boldsymbol{y}^1 \, |&\boldsymbol{y}^2 \,| & ... & |\, \boldsymbol{y}^k \end{bmatrix}.
\end{align}
We can then employ some dimension-reduction technique to extract a lower-dimensional representation of $Y$. Let $\Psi: \mathbb{R}^d \to \mathbb{R}^p$, with $p\leq d$ be a mapping defining the chosen dimension-reduction algorithm. The reduced representation of $Y$ is then given by $\overline{Y}$: 
\begin{equation}\label{reducedSnapshotMtx}
\overline{Y} = \Psi(Y) = \begin{bmatrix} \Psi(\boldsymbol{y}^1) \, |&\Psi(\boldsymbol{y}^2) \,| & ... & |& \Psi(\boldsymbol{y}^k) \end{bmatrix}\end{equation}
Different algorithms can be used to define $\Psi$, including standard Principal Component Analysis (PCA) and variants such as kernel PCA \cite{scholkopf1997kernel}, among others. The algorithm discussed in the present work is not restricted to any particular technique. For simplicity, we will assume a standard PCA, which we define using the Singular Value Decomposition (SVD). In this case, we assume \eqref{snapshotMtx} to have rows with zero mean and normalized variance.
\par Taking the SVD of $Y=USV^T$, we consider only the first $p$ singular values/vectors. Denoting the $i$-th column of $U$ as $U_i$, the first $p$ \textit{left} singular vectors $U_{1:p}$ are the \textit{principal components} of $Y$, and define a set of orthonormal axes in the reduced space whose  directions capture the majority of the variance observed in $Y$. We can project $Y$ onto the reduced space as:
\begin{equation}\label{surrogateReducedModel}
    \overline{Y} = U_{1:p}^T Y.
\end{equation}
$\overline{Y}$ is a $p \times k $ matrix whose $j$-th column expresses the $j$-th column of $Y$ in the coordinate system defined by $U_{1:p}$. 
Let us then denote by $\tildey: [0,1]^m \rightarrow \mathbb{R}^p$ the mapping that associates to $\boldtheta{}$ an output vector in the reduced space of dimension $p$, and such that $\tildey(\boldtheta{j}) = \overline{Y}_j = U_{1:p}^T\boldsymbol{y}^j$ for $j=1,\ldots,k$. We can then introduce a surrogate model $\tildePhi: [0,\,1]^m \to \mathbb{R}^p $ of $\boldsymbol{y}$ over the reduced output space as follows. We consider the following modal expansion of $\tildey$:
\begin{align}\label{reducedInterpretation}
\tildePhi(\boldtheta{}) = \sum_{i=1}^s\braket{\tildey(\boldtheta{}) | \mathcal{L}_i (\boldtheta{})}_{\rho}  \mathcal{L}_i \left( \boldtheta{}\right),
\end{align}
where $\braket{\,\cdot\, |\, \cdot \,}_{\rho}$ denotes the $\rho$-inner product over $[0,\,1]^m$ and 
 $\{\mathcal{L}_i\}_{i=1}^s$ are a set of $\rho$-orthogonal polynomials. This is the so-called \textit{generalized Polynomial Chaos Expansion (gPCE)} \cite{sudret2008global} truncated over a finite number of terms. 
 Since in the present work, we assume that $\theta_i$ are uniformly-distributed random variables, the $\mathcal{L}_i$ are Legendre polynomials \cite{sudret2008global}. 
In practice, the computation of gPCE coefficients is done using the Sparse Grids Matlab Kit \cite{piazzola.tamellini:SGK,sgmk:github}, which requires only the reduced output data $\overline{Y}$ obtained by the PCA as described above, using a sparse-grid collocation strategy to select the $\boldtheta{j}$. We refer the reader to \cite{piazzola.tamellini:SGK} for more details. 
Note that defining $\tildePhi(\boldtheta{})$ in this way allows for the fast calculation of its partial derivatives in $\theta_i$ via recurrence formulas \cite{dunster2010legendre}.

\subsection{Reduction of the input space}
 \par Observe that $\tildePhi(\boldtheta{})$ defines a map from the $m$-dimensional \textit{full input} space to the $p$-dimensional \text{reduced output space}. We now seek to obtain a corresponding $p$-dimensional input space $\widetilde{\Gamma}$. Let $\bolde{j},\,j=1,\,2,\,...,\,p$ denote the $p$-dimensional standard basis vectors; in the coordinate system defined by $U_{1:p}$, $\bolde{j}$ corresponds to the direction spanned by $U_j$. \par We aim to reduce the dimension of the input space by finding directions in $\boldsymbol{\tau}_j \in  \mathbb{R}^m$ such that a perturbation of $\boldtheta{}$ along $\boldsymbol{\tau}_j$ results in a perturbation of $\boldsymbol{y}(\boldtheta{})$ along $U_j$. Assuming such a $\boldsymbol{\tau}_j$ exists, fixing $\boldtheta{0}$ in the interior of $[0,\,1]^{m}$: 
 \begin{equation}
\boldsymbol{y}(\boldtheta{0}+\delta \boldsymbol{\tau}_j ) = \boldsymbol{y}(\boldtheta{0})+ \varepsilon U_j,\end{equation}
for $\varepsilon>0,\,\delta>0$. Observe that:
\begin{equation}\label{intuitiveArgument}
U_{1:p}^T\boldsymbol{y}(\boldtheta{0}+\delta \boldsymbol{\tau}_j ) = U_{1:p}^T\boldsymbol{y}(\boldtheta{0})+ \varepsilon U_{1:p}^T U_j \to \widetilde{\boldsymbol{y}}(\boldtheta{0}+\delta \boldsymbol{\tau}_j) -\widetilde{\boldsymbol{y}}(\boldtheta{0}) \approx \varepsilon \bolde{j}.
 \end{equation}
Replacing $\widetilde{\boldsymbol{y}}$ with its gPCE \eqref{reducedInterpretation} and exploiting its differentiability, the limit:
\begin{equation}\label{direcDerivLim}
    \lim_{\delta \to 0} \dfrac{ \tildePhi(\boldtheta{0} + \delta \boldsymbol{\tau}_j) - \tildePhi(\boldtheta{0}) }{\delta} =   \nabla \tildePhi\left(\boldtheta{0}\right) \boldsymbol{\tau}_j
\end{equation}
exists, which, together with \eqref{intuitiveArgument}, suggests:
\begin{equation}\label{reducedInputLim}
    \nabla \tildePhi\left(\boldtheta{0}\right) \boldsymbol{\tau}_j \approx \bolde{j}.
\end{equation}
Recall that:
\begin{equation}
\nabla \tildePhi(\boldtheta{0}) = \begin{bmatrix} \partial\widetilde{\Phi}_1/\partial\theta_1 & \partial\widetilde{\Phi}_1/\partial\theta_2 & ..& \partial\widetilde{\Phi}_1/\partial\theta_m  \\ 
 \partial\widetilde{\Phi}_2/\partial\theta_1 & \partial\widetilde{\Phi}_2/\partial\theta_2 & ..& \partial\widetilde{\Phi}_2/\partial\theta_m \\
\vdots & \vdots & \ddots & \vdots \\
 \partial\widetilde{\Phi}_p/\partial\theta_1 & \partial\widetilde{\Phi}_p/\partial\theta_2 & ..& \partial\widetilde{\Phi}_p/\partial\theta_m 
\end{bmatrix},
\end{equation}
with the dependence of each $\widetilde{\Phi}_i$ on $\boldtheta{0}$ understood. From the above, we observe that $\nabla \tildePhi(\boldtheta{0})$ is a $p \times m$ matrix, and $p \neq m$ in general. Therefore, $\nabla \tildePhi(\boldtheta{0})^{-1}$ will not typically exist, and, since in most cases we expect $m>p$, the system \eqref{reducedInputLim} will have infinitely many solutions in general. 
\par As an alternative, let $U_{\Phi} S_{\Phi}V_{\Phi}^T$ denote the SVD factorization of $\nabla \tildePhi(\boldtheta{0})$ such that \newline $\nabla \tildePhi(\boldtheta{0})=U_{\Phi} S_{\Phi}V_{\Phi}^T$ and define: 
\begin{equation}\label{MoorePenrose}
    \nabla\tildePhi(\boldtheta{0})^{\dag} = V_{\Phi} \begin{pmatrix} (S_{\Phi,\,1:p\times1:p})^{-1} \\ 0 \end{pmatrix} U_{\Phi}^T.
\end{equation}
The expression \eqref{MoorePenrose} is the \textit{Moore-Penrose pseudoinverse} of $\nabla \tildePhi(\boldtheta{})$\cite{golub2013matrix}. 
\newline \textbf{Definition: Reduced-order input space. } We define the $p$-dimensional \textit{reduced-order input space} $\widetilde{\Gamma}_p$ as:
\begin{equation}\label{reducedInputDirection}
\widetilde{\Gamma}_p = \text{span}\{\boldsymbol{\tau}_j\}_{j=1}^p,\, \text{ where: }\,    \boldsymbol{\tau}_j = \nabla\tildePhi(\boldtheta{0})^{\dag} \bolde{j}.
\end{equation}
\subsubsection{Moore-Penrose pseudoinverse and the need for a transformed parameter space} Let $Z_j=\{\boldsymbol{z}\in \mathbb{R}^m \,\,|\,\, \nabla\tildePhi(\boldtheta{0}) \boldsymbol{z} = \bolde{j}\}.$ An important consequence of defining the reduced-input space as \eqref{reducedInputDirection} is that \cite{golub2013matrix}:
\begin{equation}\label{minimumNormRelation}
    \| \boldsymbol{\tau}_j \|_{2} \leq \| \boldsymbol{z} \|_2 \,\,\,\, \forall \boldsymbol{z} \in Z_j.
\end{equation}
\noindent Although the system \eqref{reducedInputLim} has infinitely many solutions (the set $Z_j$) in general, $\boldsymbol{\tau}_j$ gives the \textit{minimum-norm solution}. This ensures the reduced-order input directions $\boldsymbol{\tau}_j$ are well-defined. 
\par The use of the Moore-Penrose pseudoinverse is the motivation for considering the transformed parameter vector $\boldtheta{}=\boldsymbol{\varphi}(\widehat{\boldtheta{}})$, rather than $\widehat{\boldtheta{}}$ directly. Recall that the $\widehat{\theta}_i$ refer to physical parameters, defined over intervals $\Gamma_i$, which have different numeric ranges in general. Thus, the  notion of a `minimum-norm' solution in $\Gamma$ may not be physically meaningful. By considering the transformed parameter vector $\boldtheta{}$, we ensure that each $\theta_i$ is defined over a uniform interval, and consequently, that the uniqueness-condition \eqref{minimumNormRelation} makes sense.

\section{Numerical study of HIV-prevention interventions}
To demonstrate our introduced methodology, we consider a series of HIV prevention interventions intended to prevent HIV acquisition and transmission consisting of expanding HIV testing, increasing the uptake and adherence of pre-exposure prophylaxis (PrEP) among persons without HIV, and increasing the uptake and adherence of antiretroviral therapy (ART) among persons with HIV (PWH). Since the burden of both HIV, as well associated social and structural inequities, vary by population, prevention efforts may also need to vary by population. We consider several distinct high priority transmission groups: men who have sex with men (MSM), heterosexual females (HETF), heterosexual males (HETM), and persons who inject drugs (PWID). We use HOPE, a compartmental model of HIV transmission in the United States, to illustrate the methods' ability to reduce the dimension of the input parameter space. Furthermore, we show the reduced input space can be used to rapidly develop and assess different intervention strategies.
\subsection{HOPE Model}
\par HOPE is a dynamic compartmental model that simulates the sexually active U.S. population. HOPE model includes 273 compartments corresponding to movement from the population without HIV to persons with HIV, across various disease progression states and stages of HIV care. Each model compartment is stratified by age, sex, race/ethnicity, risk level, transmission group and circumcision status. For inputs for which data were limited or uncertain, we obtained input distributions (assumed uniform) through a calibration process by varying those values within expected bounds. These inputs were calibrated so that key model outputs matched the most recent data available at the time of the analysis. Further details can be found in the technical report of the latest publication \cite{viguerie2023assessing}.
\subsection {Reduced-order inputs/outputs and interpretation}
We consider as model inputs $\boldtheta{}$, consisting of $m=10$ intervention strategies corresponding to increased levels of HIV-related prevention and care services reaching different transmission groups. Our model outputs $\boldsymbol{y}$ consist of $d=6$ HIV-related outcomes. Together, $\boldtheta{}$ and $\boldsymbol{y}$ are defined as:
\begin{equation}\label{inputsOutputs}\footnotesize
\boldtheta{}=\begin{pmatrix}
\text{PrEP (MSM)} \\ \text{PrEP (HETF)} \\ \text{PrEP (HETM)} \\ \text{PrEP (PWID)} \\ \text{ART (MSM)} \\ \text{ART (HET)} \\ \text{ART (PWID)} \\ \text{Testing (MSM)} \\ \text{Testing (HET)} \\ \text{Testing (PWID)}
\end{pmatrix},\,\,
\boldsymbol{y} = \begin{pmatrix} \text{New infections (MSM)}\\  \text{New infections (HETF)} \\  \text{New infections (HETM)} \\  \text{New infections (PWID)} \\  \text{New infections (total)}\\  \text{Total HIV-related spending (billions of 2022 USD)} \end{pmatrix}.
\end{equation}
The intervention strategies are defined as increases compared to a baseline delivery level by transmission group, beginning in the model in year 2023 and simulated through year-end 2030. The outputs refer to the cumulative totals of each outcome quantity over 2023-30. At baseline intervention levels, over the period 2023-30, our key outputs of interest are simulated as:
\begin{equation}\label{baseline}\footnotesize
\boldsymbol{y}_{\text{Baseline}} = \begin{pmatrix}
200267 & \text{New infections (MSM)}\\ 
28590  & \text{New infections (HETF)} \\ 
15689 & \text{New infections (HETM)} \\ 
16167  & \text{New infections (PWID)}   \\ 
261000 &\text{New infections (total)} \\ 
367  & \text{Total HIV-related spending (billions of 2022 USD)}\end{pmatrix}.
\end{equation}
\par Using the Sparse Grids Matlab Kit \cite{piazzola.tamellini:SGK,sgmk:github}, we generate a total of 221 points $\boldtheta{j}$, corresponding to different levels of the ten interventions \eqref{inputsOutputs}. We collect each corresponding $\boldsymbol{y}^j$ into the 6$\times$221 matrix $Y$, see \eqref{snapshotMtx}. From the SVD of $Y$, we found that considering only the first $p=4$ singular vectors/values retained over 95\% of the variance in $Y$. We then construct the map $\tildePhi$ from the intervention space to the reduced, 4-dimensional output space, following \eqref{reducedInterpretation}.
\par
\begin{figure}[t]
\centering
\includegraphics[width=10cm]{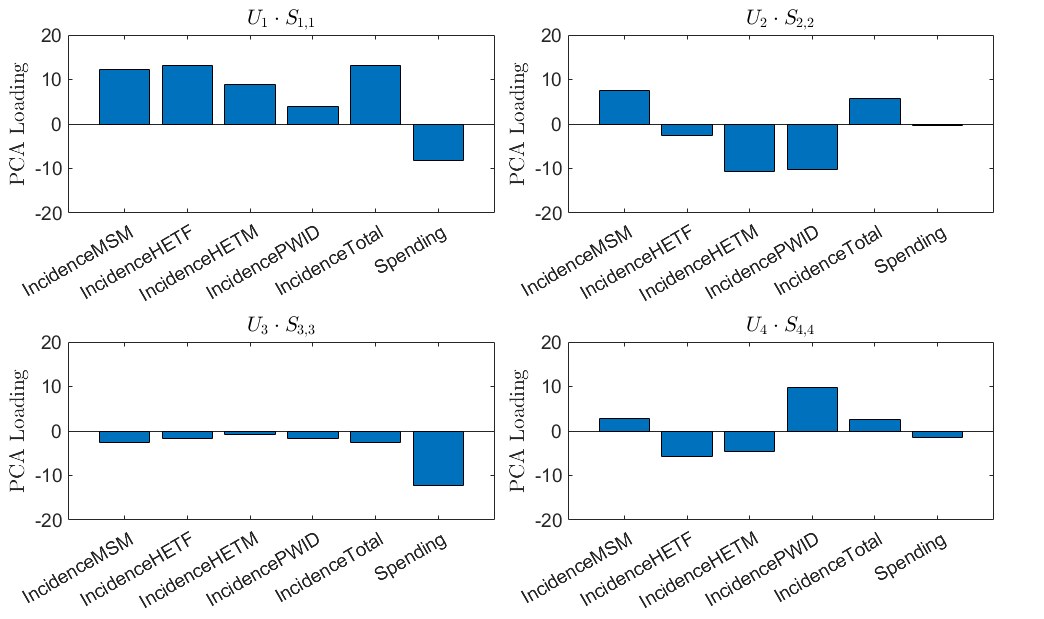}
\caption{The four identified reduced-order outputs for the HIV prevention test case.}\label{fig:pcTruncation}
\end{figure}

 The \textit{principal component loadings} for each of the 4 reduced-order output directions $U_{1:4}$ are plotted in Fig. \ref{fig:pcTruncation}. These are the entries of each vector $U_i$, weighted by the corresponding singular value $S_{i,i}$; each entry indicates how much information of each physical variable is contained in the principal component. The first principal component is strongly associated with incidence in all groups except PWID. The second component is most associated with incidence among the HETM and PWID groups. The third principal component only reflects variations in spending, while the fourth component is predominantly associated with incidence among PWID.
\par
\begin{figure}[t]
\centering
\includegraphics[width=10cm]{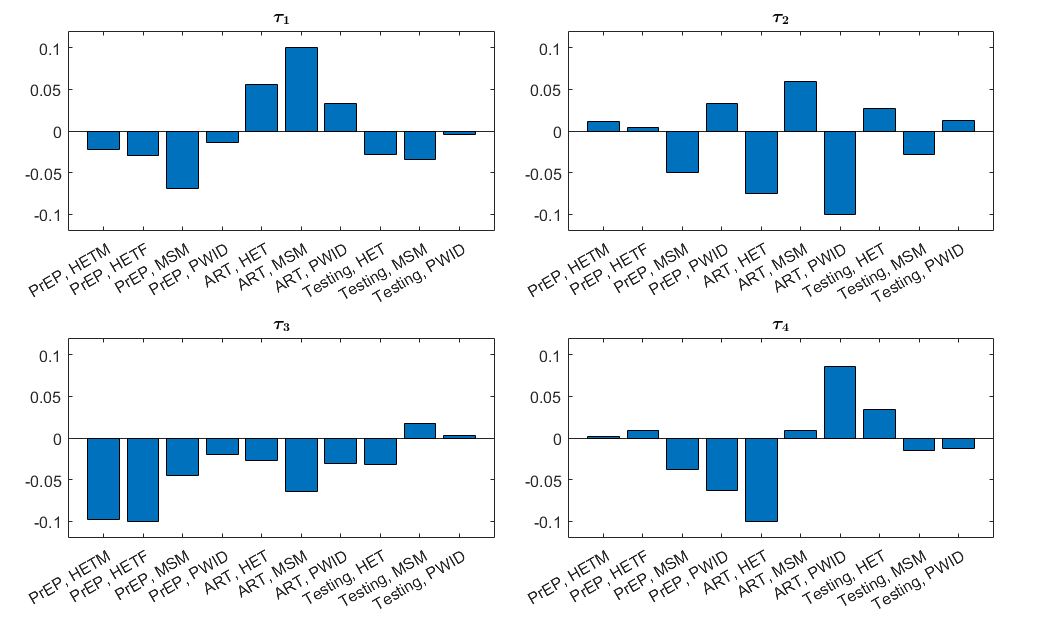}
\caption{Reduced-order inputs corresponding to each reduced-order output. }\label{fig:reducedOrderInput}
\end{figure}

 We then reduce the dimension of the input space from 10 to 4, following \eqref{reducedInputDirection} (note that $\boldtheta{0}= \dfrac{1}{2}\vec{\boldsymbol{1}}$ for convenience). The entries of each reduced-order input direction $\boldsymbol{\tau}_j$ are plotted in Fig. \ref{fig:reducedOrderInput}. The first reduced-order input direction $\boldsymbol{\tau}_1$ is dominated by PrEP and ART for MSM. The second and fourth reduced-order input directions $\boldsymbol{\tau}_{2,4}$ are primarily driven by ART among HETM and PWID ($\boldsymbol{\tau}_{2}$) and HETF and PWID ($\boldsymbol{\tau}_{4}$), respectively. Finally, the third reduced-order input direction $\boldsymbol{\tau}_3$ corresponds primarily to PrEP use among heterosexuals.
\par To confirm that the reduced-order input directions calculated by \eqref{reducedInputDirection} are correct, we compute $\boldsymbol{y}_{\boldsymbol{\tau}_j}=\boldsymbol{y}\left(\left(\boldsymbol{x}\circ \varphi^{-1}\right)\left(\boldsymbol{\theta}^0+\boldsymbol{\tau}_j \right) \right)$ for each $\boldsymbol{\tau}_j$. If the calculated $\boldsymbol{\tau}_j$ are correct, we expect that (denoting  $\boldsymbol{y}_{\boldtheta{0}}=\boldsymbol{y}\left(\left(\boldsymbol{x}\circ \varphi^{-1}\right)\left(\boldsymbol{\theta}^0\right) \right)$):
\begin{equation}\label{verificationCondition}
U_{1:4}^T\left(\boldsymbol{y}_{\boldsymbol{\tau}_j} - \boldsymbol{y}_{\boldtheta{0}} \right) \approx \bolde{j}.
\end{equation}
 Since we used a truncated SVD and a PCE, we do not expect the relation \eqref{verificationCondition} to hold exactly, but approximately. In fact, we find:
\begin{eqnarray}\begin{split}\footnotesize    
U_{1:4}^T\left(\boldsymbol{y}_{\boldsymbol{\tau}_1} - \boldsymbol{y}_{\boldtheta{0}} \right) = \begin{pmatrix}
   1.0073 \\ 0.0048 \\ 0.008 \\ 0.0024 
\end{pmatrix},\, 
U_{1:4}^T\left(\boldsymbol{y}_{\boldsymbol{\tau}_2} - \boldsymbol{y}_{\boldtheta{0}} \right) = \begin{pmatrix}
   0.0031 \\ 1.0146 \\ 0.0022 \\ 0.0004 
\end{pmatrix}, \\ \footnotesize
U_{1:4}^T\left(\boldsymbol{y}_{\boldsymbol{\tau}_3} - \boldsymbol{y}_{\boldtheta{0}} \right) = \begin{pmatrix}
   0.0007 \\ 0.0011 \\ 1.0088 \\ 0.0002 
\end{pmatrix},\,
U_{1:4}^T\left(\boldsymbol{y}_{\boldsymbol{\tau}_4} - \boldsymbol{y}_{\boldtheta{0}} \right) = \begin{pmatrix}
  -0.0001 \\ 0.0176 \\ 0.0035 \\ 1.0076 
\end{pmatrix}.
\end{split}\end{eqnarray}
\subsection{Direct use in intervention planning}
\begin{figure}[t]
\centering
\includegraphics[width=10cm]{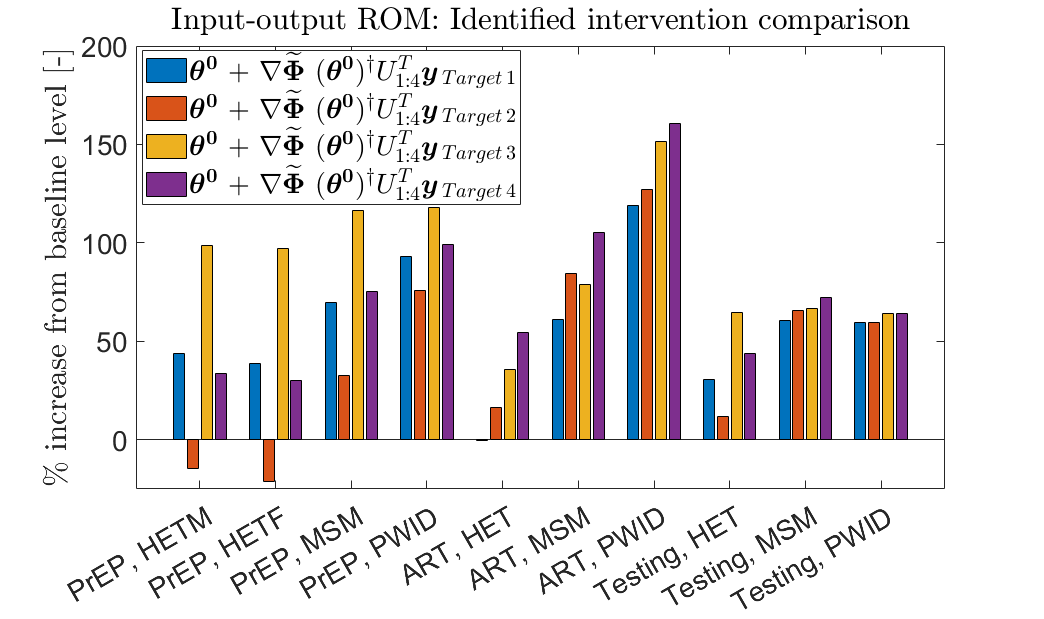}
\caption{Distinct interventions identified through input-output ROM procedure.}\label{fig:ioRomInterventions}
\end{figure}

We now apply the method to a hypothetical problem of prioritizing intervention strategies. 
We wish to find combinations of intervention strategies that meet the following hypothetical output targets (targets on HIV-related outcomes, units as in \eqref{inputsOutputs}):
\begin{equation}\label{yTargets}\footnotesize
\boldsymbol{y}_{\text{Target1}} = \begin{pmatrix}150000 \\ 24000 \\ 14000  \\ 12000  \\ 200000 \\ 425 \end{pmatrix},\boldsymbol{y}_{\text{Target2}} = \begin{pmatrix}150000\\ 24000 \\ 14000 \\ 12000 \\ 200000\\ 380  \end{pmatrix},
\boldsymbol{y}_{\text{Target3}} = \begin{pmatrix}130000\\ 22000 \\ 12000 \\ 11000 \\ 175000\\ 475  \end{pmatrix},
\boldsymbol{y}_{\text{Target4}} = \begin{pmatrix}130000\\ 22000 \\ 12000 \\ 11000 \\ 175000\\ 425  \end{pmatrix}.
\end{equation}
We can identify distinct intervention combinations that achieve these targets through with $\nabla \tildePhi{}(\boldtheta{0})^{\dag}$. Note $\boldsymbol{y}_{\text{Target} j}$ can be expressed approximately as a linear combination of the $U_i$. From \eqref{reducedInputDirection}:
\begin{equation}
\boldsymbol{y}_{\text{Target} j} \approx \sum_{i=1}^4 a_i U_i \to U_{1:4}^T\boldsymbol{y}_{\text{Target} j} \approx \sum_{i=1}^4 a_i \boldsymbol{e}_i  \to \nabla \tildePhi(\boldtheta{0})^{\dag}U_{1:4}^T\boldsymbol{y}_{\text{Target} j} \approx \sum_{i=1}^4 a_i \boldsymbol{\tau}_i.
\end{equation}
Evaluating:
\begin{eqnarray}
    \boldsymbol{y}_{\text{Eval }j}=\boldsymbol{y}\left(\left(\boldsymbol{x}\circ \varphi^{-1}\right)\left(\boldsymbol{\theta}^0+ \nabla \tildePhi(\boldtheta{0} )^{\dag}U_{1:4}^T\boldsymbol{y}_{\text{Target }j}   \right)\right) = \boldsymbol{y}\left(\left(\boldsymbol{x}\circ \varphi^{-1}\right)\left(\boldsymbol{\theta}^0+ \bigSigma_{i=1}^4 a_i \boldsymbol{\tau}_i   \right)\right)
\end{eqnarray}
for each $\boldsymbol{y}_{\text{Target } j}$ \eqref{yTargets} gives:
\begin{equation}\footnotesize
\boldsymbol{y}_{\text{Eval 1}} = \begin{pmatrix}147039\\ 24733 \\ 13447 \\ 12037 \\ 197257\\ 425  \end{pmatrix}, \boldsymbol{y}_{\text{Eval 2}} = \begin{pmatrix}145933\\ 25105 \\ 13312 \\ 11918 \\ 196629\\ 386  \end{pmatrix},
 \boldsymbol{y}_{\text{Eval 3}} = \begin{pmatrix}133822\\ 22062 \\ 11801 \\ 10659 \\ 178343\\ 472  \end{pmatrix},
 \boldsymbol{y}_{\text{Eval 4}} = \begin{pmatrix}128827\\ 21780 \\ 11682 \\ 10454 \\ 172742\\ 425  \end{pmatrix}.
\end{equation}
All intervention strategy combinations come close to achieving the targets \eqref{yTargets}, and are plotted in Fig. \ref{fig:ioRomInterventions}. Note that $\boldsymbol{y}_{\text{Target 1},\,\text{Target 2}}$ have the same HIV outcome targets, but different HIV-related spending targets. For $\boldsymbol{y}_{\text{Target 1}}$, PrEP is significantly increased from baseline levels. However, to reach the same incidence targets with a lower spending target $\boldsymbol{y}_{\text{Target 2}}$, de-emphasizes PrEP and allocates more resources toward ART. Similarly, $\boldsymbol{y}_{\text{Target 3},\,\text{Target 4}}$ have identical HIV outcome targets and different spending targets. In $\boldsymbol{y}_{\text{Target 3}}$, PrEP is heavily emphasized; $
\boldsymbol{y}_{\text{Target 4}}$ allocates comparatively fewer resources towards PrEP, and more towards ART. Nevertheless, PrEP remains an important component of $\boldsymbol{y}_{\text{Target 4}}$. We note that these results are qualitatively similar to those found in other allocation analyses  \cite{sansom2021optimal,hamilton2023achieving}. We emphasize that these are hypothetical intervention targets used for illustrative purposes to demonstrate the introduced methodology, and should not be interpreted as specific program or policy recommendations. 
\par This example demonstrates the potential utility of input/output reduced-order modeling for prioritizing intervention strategies. Different intervention combinations were identified with minimal computational cost, and did not require any sort of iterative and/or optimization procedure. Despite the low-dimension of the reduced input space, the intervention combinations identified by the input-output reduced-order model are able to account for important considerations, including identifying significantly different intervention strategies in response to changes in target spending levels. 

\section{Conclusions}

We have introduced a novel method for reducing the dimension of the input space for models with high input dimensionality. The technique works through sampling the model over a set of sparse-grid points, aggregating the model outputs of interest for each sampled point, and performing a dimension reduction on the aggregated outputs. We then build a surrogate model over the reduced output space, which maps the input parameters to the reduced outputs. Subsequently, we reduce the dimension of the input space via the construction of a mapping which associates each axis in the reduced output space with a corresponding vector in the input space. We used a national-level model of HIV transmission to demonstrate the utility of the proposed method, and found that the approach could be applied directly to construct intervention combination strategies based on hypothetical output targets, including targets related to both HIV outcomes and HIV-related spending. Furthermore, despite its reduced dimension, sophisticated dynamics, including significant changes to intervention combination strategies in response to changes in target spending levels, were captured by the reduced-order input-output space.

\par The present work is intended as a proof-of-concept, and much future work remains to develop the introduced methodology. For simplicity, we used a standard PCA to reduce the order of the output space. However, other dimension-reduction techniques exist, and may be more appropriate in certain situations. In principle, the proposed method should work with any such technique, however, further investigation is needed. To emphasize the intuition behind the approach, the method was not introduced in a mathematically rigorous manner; a more rigorous treatment, establishing existence conditions and other important properties, is necessary. Furthermore, while we emphasized public health applications herein, due to the ubiquity of models with high-dimensional input spaces in the field, the introduced technique is general in nature, and can be applied to any parameterized model satisfying very basic continuity conditions. Using the proposed method in other application areas is also an important next step. 

\medskip

\noindent\textbf{ACKNOWLEDGMENTS:} Funding for this work was provided by the Centers for Disease Control and Prevention. C. Piazzola gratefully acknowledges the support of the Alexander von Humboldt Foundation.

% \begin{thebibliography}{99}
% \bibitem{Zienkiewicz}  Zienkiewicz, O.C. and  Taylor, R.L. \textit{The finite element method}. McGraw Hill,
% Vol. I., (1989), Vol. II., (1991).
% \bibitem{Idelsohn} Idelsohn, S.R. and O\~{n}ate, E. Finite element and finite volumes. Two good friends.
% \textit{Int. J. Num. Meth. Engng.} (1994) \textbf{37}:3323--3341.
% \end{thebibliography}

%Bibliography
\bibliographystyle{unsrt}  
\bibliography{references}

\begin{thebibliography}{10}

\bibitem{luo2018proper}
Z.~Luo and G.~Chen.
\newblock {\em Proper orthogonal decomposition methods for partial differential equations}.
\newblock Academic Press, 2018.

\bibitem{quarteroni2014reduced}
A~Quarteroni, G~Rozza, et~al.
\newblock {\em Reduced order methods for modeling and computational reduction}, volume~9.
\newblock Springer, 2014.

\bibitem{jakeman2022adaptive}
J.D. Jakeman, S.~Friedman, M.S. Eldred, L.~Tamellini, A.A. Gorodetsky, and D.~Allaire.
\newblock {Adaptive experimental design for multi-fidelity surrogate modeling of multi-disciplinary systems}.
\newblock {\em Int J Num Methods Eng}, 123(12):2760--2790, 2022.

\bibitem{baddoo2023physics}
P.J. Baddoo, B.~Herrmann, B.J. McKeon, J.~Nathan~Kutz, and S.L. Brunton.
\newblock Physics-informed dynamic mode decomposition.
\newblock {\em P Roy Soc A-Math Phy}, 479(2271):20220576, 2023.

\bibitem{hess2023data}
M.W. Hess, A.~Quaini, and G.~Rozza.
\newblock {A data-driven surrogate modeling approach for time-dependent incompressible Navier-Stokes equations with dynamic mode decomposition and manifold interpolation}.
\newblock {\em Adv Comput Math}, 49(2):22, 2023.

\bibitem{viguerie2022coupled}
A.~Viguerie, G.F. Barros, M.~Grave, A.~Reali, and A.L.G.A. Coutinho.
\newblock {Coupled and uncoupled dynamic mode decomposition in multi-compartmental systems with applications to epidemiological and additive manufacturing problems}.
\newblock {\em Comput Method Appl M}, 391:114600, 2022.

\bibitem{andreuzzi2023dynamic}
F.~Andreuzzi, N.~Demo, and G.~Rozza.
\newblock {A dynamic mode decomposition extension for the forecasting of parametric dynamical systems}.
\newblock {\em SIAM J Appl Dyn Syst}, 22(3):2432--2458, 2023.

\bibitem{kim2024dimensionality}
J.~Kim, S.~Yi, and Z.~Wang.
\newblock Dimensionality reduction can be used as a surrogate model for high-dimensional forward uncertainty quantification.
\newblock {\em arXiv preprint arXiv:2402.04582}, 2024.

\bibitem{gadd2019surrogate}
C.~Gadd, W.~Xing, M.M. Nezhad, and A.A. Shah.
\newblock A surrogate modelling approach based on nonlinear dimension reduction for uncertainty quantification in groundwater flow models.
\newblock {\em Transport Porous Med}, 126:39--77, 2019.

\bibitem{ma2011kernel}
X.~Ma and N.~Zabaras.
\newblock Kernel principal component analysis for stochastic input model generation.
\newblock {\em J Comput Phys}, 230(19):7311--7331, 2011.

\bibitem{constantine2015active}
P.~G. Constantine.
\newblock {\em {Active subspaces: Emerging ideas for dimension reduction in parameter studies}}.
\newblock SIAM, 2015.

\bibitem{romor2022kernel}
F.~Romor, M.~Tezzele, A.~Lario, and G.~Rozza.
\newblock {Kernel-based active subspaces with application to computational fluid dynamics parametric problems using the discontinuous Galerkin method}.
\newblock {\em Int J Numer Meth Eng}, 123(23):6000--6027, 2022.

\bibitem{scholkopf1997kernel}
B.~Sch{\"o}lkopf, A.~Smola, and K.~M{\"u}ller.
\newblock Kernel principal component analysis.
\newblock In {\em International conference on artificial neural networks}, pages 583--588. Springer, 1997.

\bibitem{sudret2008global}
B.~Sudret.
\newblock Global sensitivity analysis using polynomial chaos expansions.
\newblock {\em Reliab Eng Syst Safe}, 93(7):964--979, 2008.

\bibitem{piazzola.tamellini:SGK}
C.~Piazzola and L.~Tamellini.
\newblock {Algorithm 1040: The Sparse Grids Matlab Kit - a Matlab implementation of sparse grids for high-dimensional function approximation and uncertainty quantification}.
\newblock {\em ACM T Math Software}, 50(1), 2024.

\bibitem{sgmk:github}
C.~Piazzola and L.~Tamellini.
\newblock {The Sparse Grids Matlab Kit - v. 23-5 Robert}, 2023.
\newblock \url{https://github.com/lorenzo-tamellini/sparse-grids-matlab-kit}.

\bibitem{dunster2010legendre}
T.M. Dunster.
\newblock Legendre and related functions, 2010.

\bibitem{golub2013matrix}
G.H. Golub and C.~F. Van~Loan.
\newblock {\em Matrix computations}.
\newblock JHU press, 2013.

\bibitem{viguerie2023assessing}
A.~Viguerie, E.U. Jacobson, K.A. Hicks, L.~Bates, J.~Carrico, A.~Honeycutt, C.~Lyles, and P.G. Farnham.
\newblock {Assessing the impact of COVID-19 on HIV Outcomes in the United States: A modeling study}.
\newblock {\em Sex Transim Dis}, pages 10--1097, 2023.

\bibitem{sansom2021optimal}
S.L. Sansom, K.A. Hicks, J.~Carrico, E.U. Jacobson, R.K. Shrestha, T.A. Green, and D.W. Purcell.
\newblock {Optimal allocation of societal HIV prevention resources to reduce HIV incidence in the United States}.
\newblock {\em Am J Public Health}, 111(1):150--158, 2021.

\bibitem{hamilton2023achieving}
D.~Hamilton, K.~Hoover, D.~Smith, K.~Delaney, L.~Wang, J.~Li, T.~Hoyte, S.~Jenness, and S.~Goodreau.
\newblock {Achieving the “Ending the HIV Epidemic in the US” incidence reduction goals among at-risk populations in the South}.
\newblock {\em BMC Public Health}, 23(1):716, 2023.

\end{thebibliography}

\end{document}